\documentclass[12pt]{article}
\usepackage{graphicx}
\usepackage{amsmath}
\usepackage[dvips]{epsfig}
\usepackage{amssymb}

\makeatletter
\renewcommand{\section}{\@startsection
  {section}%
  {2}%
  {0mm}%
  {\baselineskip}%
  {0.5 \baselineskip}%
  {\centering}}
\makeatother
\begin{document}

\title {Some new identities on the twisted $(h, q)$-Euler numbers and $q$-Bernstein  polynomials }
\author{D. V. Dolgy$^1$,  D. J. Kang$^2$,  T.  Kim$^3$, and B. Lee$^4$\\[0.5cm]
             $^1$Hanrimwon,
 Kwangwoon University, Seoul 139-701,  Korea\\\\
 $^2$ Information Technology Service,\\           
   Kyungpook National University,  Taegu 702-701, Korea \\ \\ 
             $^3$Division of General Education-Mathematics,\\
 Kwangwoon University, Seoul 139-701,  Korea\\\\
  $^4$Department of Wireless Communications Engineering,\\
  Kwangwoon University, Seoul 139-701,  Korea\\\\
   }
\date{}
\maketitle

 {\footnotesize {\bf Abstract}\hspace{1mm}
In this paper we give some interesting relationship between the
$q$-Bernstein polynomials and the twisted  $(h, q)$-Euler numbers by
using fermionic $p$-adic $q$-integral on $\Bbb Z_p$.



\bigskip
{\bf\section{Introduction }}
\bigskip

Let $p$ be a fixed odd prime number. Throughout this paper, we
always make use of the following notations: $\Bbb Z$ denotes the
ring of rational integers,  $\mathbb{Z}_p$ denotes the ring of
$p$-adic rational integers, $\mathbb{Q}_p$ denotes the field of
$p$-adic rational numbers, and $\mathbb{C}_p$ denotes the
completion of algebraic closure of $\mathbb{Q}_p$, respectively.
 Let $\Bbb N $
be the set of natural numbers and $\Bbb Z_+ = \Bbb N \cup \{0 \}.$
Let $C_{p^n}=\{ \zeta| \zeta^{p^n}=1\}$ be the cyclic group of
order
 $p^n$ and let $$ T_p= \lim_{n \rightarrow \infty} C_{p^n} =\cup_{n\geq
 1}C_{p^n}.$$
The $p$-adic absolute value is defined by $|x|_p =
\dfrac{1}{p^r}$, where $x=p^r \dfrac{s}{t}$( $ r \in \Bbb Q $ and
$s, t \in \Bbb Z$  with $(s,t)=(p,s)=(p,t)=1$). In this paper we
assume  that $q\in \Bbb C_p$ with $|q-1|_p<1$ as an indeterminate.
The $q$-number is defined by
$$[x]_q =\frac{1-q^x}{1-q}, \text{ see [1-18].}$$ Note that $\lim_{ q \rightarrow 1}[x]_q=x$.
 For
$$ f \in UD(\mathbb{Z}_p)=\{ f | f :\mathbb{Z}_p \to \mathbb{C}_p \text { is uniformly differentiable function} \}, $$
the fermionic $p$-adic $q$-integral on $\mathbb{Z}_p$  is defined
as
$$I_{-q}(f)=\int_{\mathbb{Z}_p}f(x) d\mu_{-q}(x)=
\lim_{N \to \infty} \dfrac{1+q}{1+q^{p^N}} \sum_{ x=0}^{ p^N-1}
f(x)(-q)^x, \mbox{ see [1-5] }. \eqno(1)$$
 From (1), we note that
$$ q^n I_{-q}(f_n)=  (-1)^{n}I_{-q}(f)+ [2]_q  \sum_{l=0}^{n-1} (-1)^{n-1-l}q^l f(l),
$$ where $ f_n(x)=f(x+n)$ for $ n \in \Bbb N.$
 For $k, n \in \Bbb Z_+$ and $ x \in [0, 1]$,
 Kim defined $q$-Bernstein polynomials, which are different $q$-Bernstein polynomials of
 Phillips, as follows:
$$  B_{k, n}(x, q) = \binom{n}{k}[x]_q^k [1-x]_{q^{-1}}^{n-k}, \text{ see [4].} \eqno(2) $$

In [9], the $p$-adic extension of (2) is given by
$$  B_{k, n}(x, q) = \binom{n}{k}[x]_q^k [1-x]_{q^{-1}}^{n-k},
 \text{ where } x \in \Bbb Z_p, \text{ and } n, k \in \Bbb Z_+. \eqno(3) $$

For $h \in \mathbb{Z}$ and $\zeta \in T_p$, let us consider the
twisted  $(h, q)$-Euler polynomials as follows:
$$ {E}_{n,q, \zeta}^{(h)}(x)=\int_{\mathbb{Z}_p} [x+y]_q^n  \zeta^y q^{(h-1)y} d\mu_{-q}(y), \text{ for } n \in \Bbb Z_+. \eqno(4)
$$
In the special case, $x=0$,
${E}_{n,q,\zeta}^{(h)}(0)={E}_{n,q,\zeta}^{(h)}$ are called the
$n$-th twisted $(h, q)$-Euler numbers.

In this paper  we investigate some relations between the
$q$-Bernstein polynomials and the twisted $(h, q)$-Euler numbers.
From these relations, we derive some interesting identities on the
the twisted $(h, q)$-Euler numbers and polynomials.

\bigskip
{\bf\section{On the twisted $(h, q)$-Euler numbers and  polynomials}}
\bigskip

From (4), we note that
$$ \aligned  {E}_{n,q, \zeta}^{(h)}(x)
&= \int_{\mathbb{Z}_p} [x+y]_q^n  \zeta^y q^{(h-1)y}
d\mu_{-q}(y)\\
&= \dfrac{[2]_q}{(1-q)^n}\sum_{l=0}^n \binom{n}{l}(-1)^l q^{lx}
\left( \dfrac{1}{1+q^{h+l}\zeta} \right)
\\
&= [2]_q \sum_{m=0}^\infty (-1)^m \zeta^m q^{hm}[x+m]_q^n,
\endaligned $$
and
$$ \aligned  {E}_{n,q, \zeta}^{(h)}(x)
&= \int_{\mathbb{Z}_p} [x+y]_q^n  \zeta^y q^{(h-1)y}
d\mu_{-q}(y)\\
&=  \sum_{l=0}^n \binom nl [x]_q^{n-l} q^{lx}\int_{\mathbb{Z}_p}
[y]_q^l  \zeta^y q^{(h-1)y} d\mu_{-q}(y)
\\
&= \sum_{l=0}^n \binom nl [x]_q^{n-l} q^{lx}{E}_{l,q,
\zeta}^{(h)}.
\endaligned \eqno(5)$$

Therefore, we obtain the following theorem.

\bigskip

{ \bf Theorem 1.}  For $ n \in \Bbb Z_+$ and  $\zeta \in T_p$, we
have
$$ \aligned  {E}_{n,q, \zeta}^{(h)}(x)
&= [2]_q \sum_{m=0}^\infty (-1)^m \zeta^m q^{hm}[x+m]_q^n.
\endaligned $$
Furthermore,
$$ \aligned  {E}_{n,q,\zeta}^{(h)}(x) &= \sum_{l=0}^n \binom nl [x]_q^{n-l} q^{lx}  {E}_{l,q,\zeta}^{(h)}\\
&=([x]_q+ q^x {E}_{q,\zeta}^{(h)})^n, \endaligned $$ with  usual
convention  about replacing $ ( {E}_{q, \zeta }^{(h)})^n$ by
${E}_{n,q, \zeta}^{(h)}$.
\bigskip

Let $$F_{q, \zeta}^{(h)}(t, x)=\sum_{n=0}^\infty
{E}_{n,q,\zeta}^{(h)}(x) \dfrac{t^n}{n!}.$$  Then we see that
$$ F_{q,\zeta}^{(h)}(t,x)=[2]_q \sum_{m=0}^\infty (-1)^m \zeta^m q^{mh}e^{[x+m]_qt}. \eqno(6) $$
In the special case, $x=0$, let  $ F_{q, \zeta}^{(h)}(t, 0)=F_{q,
\zeta} ^{(h)}(t).$

By (4), we get

$$ \aligned  {E}_{n,q^{-1}, \zeta^{-1}}^{(h)}(1-x)
&= \int_{\mathbb{Z}_p} [1-x+y]_{q^{-1}}^n  \zeta^{-y} q^{-(h-1)y}
d\mu_{-q^{-1}}(y)\\
&=  \dfrac{[2]_q}{(1-q^{-1})^n}\sum_{l=0}^n \binom{n}{l}(-1)^l
q^{h-1} \zeta  \left( \dfrac{q^{lx}}{1+q^{h+l} \zeta} \right)
\\
&= (-1)^n \zeta q^{n+h-1} \left(
\dfrac{[2]_q}{(1-q)^n}\sum_{l=0}^n \binom{n}{l}(-1)^l
\dfrac{q^{lx}}{1+q^{h+l}\zeta} \right)
\\
&=  (-1)^n \zeta q^{n+h-1} {E}_{n, q, \zeta}^{(h)}(x).
\endaligned$$

Therefore, we obtain the following theorem.

{ \bf Theorem 2.}  Let $ F_{q,\zeta}^{(h)}(t,x)= \sum_{n=0}^\infty
{E}_{n,q,\zeta}^{(h)}(x) \dfrac{t^n}{n!}. $ Then we have
$$ F_{q^{-1},\zeta^{-1}}^{(h)}(t,1-x)= \zeta q^{h-1}
F_{q,\zeta}^{(h)}(-qt,x).$$ Moreover,
$$  {E}_{n,q^{-1}, \zeta^{-1}}^{(h)}(1-x)
= (-1)^n \zeta q^{n+h-1} {E}_{n, q, \zeta}^{(h)}(x) \text{ for } n
\in \Bbb Z_+ .$$
\bigskip

From (6), we note that
$$q^h \zeta  F_{q,\zeta}^{(h)}(t,1)+F_{q,\zeta}^{(h)}(t)=[2]_q. \eqno(7)$$
By (7), we get the following recurrence formula:
$${E}_{0, q, \zeta}^{(h)} =\dfrac{1+q}{1+  q^h \zeta}, \text{ and }    q^h \zeta  {E}_{n, q, \zeta}^{(h)}(1) +  {E}_{n,
q, \zeta}^{(h)}=0   \mbox{ if } n>0.
 \eqno(8)
$$
By (8) and Theorem 1, we obtain the following theorem.

\bigskip
{ \bf Theorem 3.}  For $ n \in \Bbb Z_+$ and  $\zeta \in T_p$, we
have
$${E}_{0, q, \zeta}^{(h)} =\dfrac{1+q}{1+  q^h \zeta}, \text{ and }    q^h \zeta  (q {E}_{ q, \zeta}^{(h)}+1)^n +  {E}_{n,
q, \zeta}^{(h)}=0   \mbox{ if } n>0,
$$ with  usual
convention  about replacing $ ( {E}_{q, \zeta }^{(h)})^n$ by
${E}_{n,q, \zeta}^{(h)}$.
\bigskip

From Theorem 3, we note that

$$
 \aligned & q^{2h} \zeta^2  {E}_{n,q, \zeta}^{(h)}( 2 )
 - \dfrac{1+q}{1+q^h\zeta} q^{2h}\zeta^2- \dfrac{1+q}{1+q^h \zeta} q^{h}\zeta \\
 & = q^{2h} \zeta^2  \sum_{l=0}^ n \binom nl q^l ( q {E}_{q, \zeta}^{(h)}+ 1)^l -
  \dfrac{1+q}{1+q^h\zeta} q^{2h}\zeta^2- \dfrac{1+q}{1+q^h \zeta} q^{h}\zeta\\
   & = q^{2h} \zeta^2  \sum_{l=1}^ n \binom nl q^l ( q {E}_{q, \zeta}^{(h)}+ 1)^l - \dfrac{1+q}{1+q^h \zeta} q^{h}\zeta\\
 &=   -q^{h} \zeta \sum_{l=1}^ n \binom nl q^l   {E}_{l,q, \zeta}^{(h)} - \dfrac{1+q}{1+q^h \zeta} q^{h}\zeta \\
 & =  -q^{h} \zeta  \sum_{l=0}^ n \binom nl q^l   {E}_{l,q, \zeta}^{(h)} \\
 &=- q^{h} \zeta {E}_{n,q, \zeta}^{(h)}( 1 ) = {E}_{n,q, \zeta}^{(h)} \text{ if } n>0.
 \endaligned $$

Therefore, we obtain the following theorem.

\bigskip

{ \bf Theorem 4.}  For $ n  \in \Bbb N$,  we have
$$   q^{2h} \zeta^2  {E}_{n,q, \zeta}^{(h)}( 2 )={E}_{n,q, \zeta}^{(h)} +
 \dfrac{1+q}{1+q^h\zeta} q^{2h}\zeta^2 + \dfrac{1+q}{1+q^h \zeta} q^{h}\zeta .$$

\bigskip

By Theorem 2, we see that

$$
 \aligned   & q^{h-1} \zeta
  \int_{\mathbb{Z}_p} [1-x]_{q^{-1}}^n q^{(h-1)x} \zeta^x
  d\mu_{-q}(x)\\
  &=(-1)^n q^{n+h-1} \zeta  \int_{\mathbb{Z}_p} [x-1]_{q}^n q^{(h-1)x} \zeta^x d\mu_{-q}(x) \\
 &= (-1)^n  q^{n+h-1}  \zeta {E}_{n,q, \zeta}^{(h)}( -1 )={E}_{n,q^{-1}, \zeta^{-1}}^{(h)}(2).
 \endaligned \eqno(9)
$$

Therefore,  we obtain the following theorem.

\bigskip

{ \bf Theorem 5.}  For $ n \in \Bbb Z_+$,  we have
$$ q^{h-1}  \zeta \int_{\mathbb{Z}_p} [1-x]_{q^{-1}}^n q^{(h-1)x} \zeta^x d\mu_{-q}(x)
= {E}_{n,q^{-1}, \zeta^{-1}}^{(h)}(2).$$

\bigskip

Let $ n \in \Bbb N.$ By Theorem 5 and Theorem 4, we get
$$
 \aligned   & q^{h-1} \zeta
  \int_{\mathbb{Z}_p} [1-x]_{q^{-1}}^n q^{(h-1)x} \zeta^x
  d\mu_{-q}(x)\\
  &= q^{2h} \zeta^2 {E}_{n,q^{-1}, \zeta^{-1}}^{(h)} +
q^{h-1} \zeta \left(  \dfrac{1+q}{1+q^h\zeta} \right)+ q^{2h-1}
\zeta^2 \left( \dfrac{1+q}{1+q^h \zeta} \right)  .
 \endaligned \eqno(10)
$$

From (10), we have

$$
  \int_{\mathbb{Z}_p} [1-x]_{q^{-1}}^n q^{(h-1)x} \zeta^x
  d\mu_{-q}(x)= q^{h+1} \zeta {E}_{n,q^{-1}, \zeta^{-1}}^{(h)} +
 \left(  \dfrac{1+q}{1+q^h\zeta} \right)+ q^{h} \zeta \left(
\dfrac{1+q}{1+q^h \zeta} \right) .
$$

Therefore, we obtain the following corollary.

\bigskip

{ \bf Corollary 6.}  For $ n \in \Bbb N$,  we have $$
  \int_{\mathbb{Z}_p} [1-x]_{q^{-1}}^n q^{(h-1)x} \zeta^x
  d\mu_{-q}(x)= q^{h+1} \zeta {E}_{n,q^{-1}, \zeta^{-1}}^{(h)} + [2]_q .
$$

\bigskip

For  $x \in \Bbb Z_p$, the $p$-adic analogues of  $q$-Bernstein
polynomials are given  by
$$B_{k,n}(x, q)= {\binom{n}{k}} [x]_q^k [1-x]_{q^{-1}}^{n-k},  \text{
where } n, k \in \Bbb Z_+. \eqno(11)
$$
By (11), we get the symmetry of $q$-Bernstein polynomials as
follows:
$$ B_{k,n}(x, q)=B_{n-k, n}(1-x, q^{-1}), \text{ see [7]}.
\eqno(12) $$ Thus, by Corollary 6, (11), and (12), we see that

$$  \aligned   \int_{\mathbb{Z}_p }    B_{k, n}(x,q) q^{(h-1)x}
\zeta^x d\mu_{-q}(x) &=\int_{\mathbb{Z}_p }    B_{n-k, n}(1-x,
q^{-1}) q^{(h-1)x} \zeta^x
d\mu_{-q}(x) \\
  & =  \binom nk \sum_{l=0}^{k} \binom{k}{l}(-1)^{k+l}
  \int_{\mathbb{Z}_p } [1-x]_{q^{-1}}^{n-l}  q^{(h-1)x} \zeta^x d\mu_{-q} (x) \\
  &=  \binom nk \sum_{l=0}^{k} \binom{k}{l}(-1)^{k+l} \left(  q^{h+1} \zeta {E}_{n-l,q^{-1}, \zeta^{-1}}^{(h)} + [2]_q  \right).\endaligned
$$

For $n, k \in \Bbb Z_+$ with $n>k$, we have
$$  \aligned   & \int_{\mathbb{Z}_p }    B_{k, n}(x,q) q^{(h-1)x}
\zeta^x
d\mu_{-q}(x)\\
  & = \binom nk \sum_{l=0}^{k} \binom{k}{l}(-1)^{k+l} \left(  q^{h+1} \zeta {E}_{n,q^{-1}, \zeta^{-1}}^{(h)} + [2]_q  \right)\\
  &=\left \{\begin{array}{ll}
   q^{h+1} \zeta {E}_{n,q^{-1}, \zeta^{-1}}^{(h)} + [2]_q , & \mbox{ if } k=0, \\ \\
  q^{h+1} \zeta \binom{n}{k} \sum_{l=0}^{k}
\binom{k}{l}(-1)^{k+l}  q^{h+1} \zeta {E}_{n-l,q^{-1},
\zeta^{-1}}^{(h)}, & \mbox{ if } k
> 0.
\end{array} \right.
 \endaligned \eqno(13)
$$
Let us take the fermionic $q$-integral on $\Bbb Z_p$ for the
$q$-Bernstein polynomials of degree $n$ as follows:

$$  \aligned   \int_{\mathbb{Z}_p }    B_{k, n}(x,q) q^{(h-1)x} \zeta^x
d\mu_{-q}(x)
  & = \binom nk  \int_{\mathbb{Z}_p }   [x]_q^k [1-x]_{q^{-1}}^{n-k}
  q^{(h-1)x}\zeta^x d\mu_{-q}(x)\\
  &= \binom{n}{k} \sum_{l=0}^{n-k} \binom{n-k}{l}(-1)^{l}
{E}_{l+k,q, \zeta}^{(h)}.
 \endaligned \eqno(14)
$$

Therefore, by (13) and (14),  we obtain the following theorem.

\bigskip

{ \bf Theorem 7.}  Let $ n, k \in \Bbb Z_+$ with $n>k$. Then   we
have
$$  \aligned   & \int_{\mathbb{Z}_p }    B_{k, n}(x,q) q^{(h-1)x}\zeta^x
d\mu_{-q}(x) \\
 &=\left \{\begin{array}{ll}
   q^{h+1} \zeta {E}_{n,q^{-1}, \zeta^{-1}}^{(h)} + [2]_q , & \mbox{ if } k=0, \\ \\
 q^{h+1}
\zeta  \binom{n}{k} \sum_{l=0}^{k} \binom{k}{l}(-1)^{k+l}
{E}_{n-l,q^{-1}, \zeta^{-1}}^{(h)}, & \mbox{ if } k
> 0.
\end{array} \right.
 \endaligned
$$
Moreover,
$$  \aligned  &   \sum_{l=0}^{n-k} \binom{n-k}{l}(-1)^{l}
{E}_{l+k,q, \zeta}^{(h)} \\
&=\left \{\begin{array}{ll}
   q^{h+1} \zeta {E}_{n,q^{-1}, \zeta^{-1}}^{(h)} + [2]_q , & \mbox{ if } k=0, \\ \\
 q^{h+1} \zeta \sum_{l=0}^{k} \binom{k}{l}(-1)^{k+l} {E}_{n-l,q^{-1}, \zeta^{-1}}^{(h)},
& \mbox{ if } k
> 0.
\end{array} \right.
 \endaligned
$$
\bigskip

 Let $n_1, n_2, k \in \Bbb Z_+$ with $n_1+n_2 > 2k.$  Then we get
$$  \aligned   & \int_{\mathbb{Z}_p } B_{k, n_1}(x,q) B_{k, n_2}(x,q) q^{(h-1)x} \zeta^x  d\mu_{-q}
(x)\\
  & =   \binom {n_1}{k} \binom {n_2}{k}\sum_{l=0}^{2k} \binom{2k}{l}(-1)^{l+2k}
   \int_{\mathbb{Z}_p } [1-x]_{q^{-1}}^{n_1+n_2-l} q^{(h-1)x}  \zeta^x  d\mu_{-q} (x)  \\
  &=  \binom {n_1}{k} \binom {n_2}{k}\sum_{l=0}^{2k} \binom{2k}{l}(-1)^{l+2k}
    \left( { q^{h+1}}\zeta {E}_{n_1+n_2-l,q^{-1}, \zeta^{-1}}^{(h)}+  {[2]_q} \right) \\
  &= \left \{\begin{array}{ll}
  { q^{h+1}}\zeta {E}_{n_1+n_2,q^{-1}, \zeta^{-1} }^{(h)}+   [2]_q, & \mbox{ if } k=0, \\ \\
\binom {n_1}{k} \binom {n_2}{k}\sum_{l=0}^{2k}
\binom{2k}{l}(-1)^{2k+l}  \left( { q^{h+1}}\zeta
{E}_{n_1+n_2-l,q^{-1}, \zeta^{-1}}^{(h)}+  {[2]_q} \right), &
\mbox{ if } k \neq 0.
\end{array} \right. \endaligned \eqno(15)
$$

Therefore,  by (15), we obtain the following theorem.

\bigskip

{ \bf Theorem 8.} For $n_1, n_2, k \in \Bbb Z_+$ with $n_1+n_2 >
2k,$    we have
$$  \aligned   & \int_{\mathbb{Z}_p } B_{k, n_1}(x,q) B_{k, n_2}(x,q) q^{(h-1)x}  \zeta^x   d\mu_{-q}
(x) \\
&= \left \{\begin{array}{ll}
  { q^{h+1}} \zeta{E}_{n_1+n_2,q^{-1}, \zeta^{-1} }^{(h)}+   [2]_q, & \mbox{ if } k=0, \\ \\
 { q^{h+1}}\zeta \binom {n_1}{k} \binom {n_2}{k}\sum_{l=0}^{2k}
\binom{2k}{l}(-1)^{2k+l} {E}_{n_1+n_2-l,q^{-1}, \zeta^{-1}}^{(h)},
& \mbox{ if } k \neq 0.
\end{array} \right. \endaligned
$$
\bigskip

From the binomial theorem, we can derive the following equation.

$$  \aligned   & \int_{\mathbb{Z}_p } B_{k, n_1}(x,q) B_{k, n_2}(x,q) q^{(h-1)x}  \zeta^x  d\mu_{-q}
(x)\\
  & =   \binom {n_1}{k} \binom {n_2}{k}\sum_{l=0}^{n_1+n_2-2k}(-1)^{l} \binom{n_1+n_2-2k}{l}
   \int_{\mathbb{Z}_p } [x]_q^{2k+l}  q^{(h-1)x} \zeta^x  d\mu_{-q} (x)  \\
  &=    \binom {n_1}{k} \binom {n_2}{k}\sum_{l=0}^{n_1+n_2-2k} (-1)^{l} \binom{n_1+n_2-2k}{l}
  E_{2k+l, q, \zeta}^{(h)}.\endaligned
  \eqno(16)
$$
Thus, by Theorem 8 and (16), we obtain the following corollary.

\bigskip

{ \bf Corollary 9.} Let $n_1, n_2, k \in \Bbb Z_+$ with $n_1+n_2 >
2k.$ Then   we have
$$  \aligned   & \sum_{l=0}^{n_1+n_2-2k} (-1)^{l} \binom{n_1+n_2-2k}{l}
  E_{2k+l, q}^{(h)} \\
  &  =  \left \{\begin{array}{ll}
  { q^{h+1}}\zeta {E}_{n_1+n_2,q^{-1}, \zeta^{-1}}^{(h)}+  [2]_q, & \mbox{ if } k=0, \\ \\
{ q^{h+1}}\zeta \sum_{l=0}^{2k} \binom{2k}{l}(-1)^{2k+l}
{E}_{n_1+n_2-l,q^{-1}, \zeta^{-1}}^{(h)}, & \mbox{ if } k > 0.
\end{array} \right. \endaligned
$$
\bigskip

For $ x \in \Bbb Z_p$ and $ s \in \Bbb N$ with $s \geq 2,$  let
$n_1, n_2, \ldots, n_s , k \in \Bbb Z_+$ with $n_1+ \cdots + n_s >
sk.$  Then we take the fermionic $p$-adic $q$-integral on $\Bbb
Z_p$ for the $q$-Bernstein polynomials of degree $n$ as follows:
$$ \aligned   & \int_{\mathbb{Z}_p }   \underset{ s-\mbox{times}}{\underbrace{
B_{k, n_1}(x, q)  \cdots B_{k, n_s}(x, q)  } }  q^{(h-1)x} \zeta^x
d\mu_{-q}(x) \\
&= \binom {n_1}{k} \cdots \binom {n_s}{k} \int_{\mathbb{Z}_p }
[x]_q^{sk} [1-x]_{q^{-1}}^{n_1+ \cdots +n_s-sk} q^{(h-1)x} \zeta^x
d\mu_{-q}(x) \\
&= \binom {n_1}{k} \cdots \binom {n_s}{k} \sum_{l=0}^{sk}
\binom{sk}{l} (-1)^{l+sk}  \int_{\mathbb{Z}_p
}[1-x]_{q^{-1}}^{n_1+ \cdots +n_s-l} q^{(h-1)x} \zeta^x d\mu_{-q}(x)\\
&=  \binom {n_1}{k} \cdots \binom {n_s}{k} \sum_{l=0}^{sk}
\binom{sk}{l} (-1)^{l+sk}  \left( q^{h+1} \zeta {E}_{n_1+
\cdots+ n_s-l,q^{-1}, \zeta^{-1}}^{(h)}+[2]_q \right) \\
&=  \left \{\begin{array}{ll}
  { q^{h+1}} \zeta {E}_{n_1+ \cdots +n_s,q^{-1},\zeta^{-1}}^{(h)}+  [2]_q, & \mbox{ if } k=0, \\ \\
  q^{h+1}  \zeta \binom {n_1}{k} \cdots \binom {n_s}{k}
\sum_{l=0}^{sk} \binom{sk}{l} (-1)^{l+sk}  {E}_{n_1+ \cdots+
n_s-l,q^{-1}, \zeta^{-1}}^{(h)}, & \mbox{ if } k > 0.
\end{array} \right. \endaligned \eqno(17)
$$

Therefore, by (17), we obtain the following theorem.

\bigskip

{ \bf Theorem 10.} For $ s \in \Bbb N$ with $s \geq 2$,   let
$n_1, n_2, \ldots, n_s , k \in \Bbb Z_+$ with $n_1+ \cdots + n_s >
sk$. Then we get
$$ \aligned   & \int_{\mathbb{Z}_p }   \underset{ s-\mbox{times}}{\underbrace{
B_{k, n_1}(x, q)  \cdots B_{k, n_s}(x, q)  } }  q^{(h-1)x} \zeta^x
d\mu_{-q}(x) \\
&= \left \{\begin{array}{ll}
  { q^{h+1}} \zeta {E}_{n_1+ \cdots +n_s,q^{-1},\zeta^{-1}}^{(h)}+  [2]_q, & \mbox{ if } k=0, \\ \\
{ q^{h+1}} \zeta  \binom {n_1}{k} \cdots \binom {n_s}{k}
\sum_{l=0}^{sk} \binom{sk}{l} (-1)^{l+sk} {E}_{n_1+ \cdots+
n_s-l,q^{-1}, \zeta^{-1}}^{(h)}, & \mbox{ if } k > 0.
\end{array} \right. \endaligned \eqno(18)
$$

 \bigskip

By the definition of $q$-Bernstein polynomials and the binomial
theorem, we easily get
$$ \aligned   & \int_{\mathbb{Z}_p }   \underset{ s-\mbox{times}}{\underbrace{
B_{k, n_1}(x, q)  \cdots B_{k, n_s}(x, q)  } }  q^{(h-1)x} \zeta^x
d\mu_{-q}(x) \\
&=  \binom {n_1}{k} \cdots \binom {n_s}{k} \sum_{l=0}^{n_1+ \cdots
+ n_s-sk}(-1)^{l} \binom{n_1+ \cdots + n_s-sk}{l}
   \int_{\mathbb{Z}_p } [x]_q^{sk+l}  q^{(h-1)x} \zeta^z d\mu_{-q} (x)  \\
  &=    \binom {n_1}{k} \cdots \binom {n_s}{k} \sum_{l=0}^{n_1+ \cdots
+ n_s-sk}(-1)^{l} \binom{n_1+ \cdots + n_s-sk}{l}
 { E}_{sk+l, q, \zeta}^{(h)}.\endaligned
 $$

Therefore,  we have the following corollary.

\bigskip

{ \bf Corollary 11.} For  $ \zeta \in T_p,  s \in \Bbb N$ with $s
\geq 2$,  let  $n_1, n_2, \ldots, n_s , k \in \Bbb Z_+$ with $n_1+
\cdots + n_s > sk.$  Then  we have
$$
\aligned   &   \sum_{l=0}^{n_1+ \cdots + n_s-sk}(-1)^{l}
\binom{n_1+ \cdots + n_s-sk}{l}
 { E}_{sk+l, q, \zeta}^{(h)} \\
  &  = \left \{\begin{array}{ll}
  { q^{h+1}} \zeta {E}_{n_1+ \cdots +n_s,q^{-1},\zeta^{-1}}^{(h)}+  [2]_q, & \mbox{ if } k=0, \\ \\
{ q^{h+1}} \zeta  \sum_{l=0}^{sk} \binom{sk}{l} (-1)^{l+sk}
{E}_{n_1+ \cdots+ n_s-l,q^{-1}, \zeta^{-1}}^{(h)}, & \mbox{ if } k
> 0.
\end{array} \right. \endaligned
$$

 \end{document}